\documentclass{amsart}
\usepackage{amsmath} 
\usepackage{amssymb}

\theoremstyle{plain}
\newtheorem{thrm}{Theorem}

\newtheorem{cject}[thrm]{Conjecture}

\theoremstyle{definition}
\newtheorem{dfin}[thrm]{Definition}
\newtheorem{prlem}[thrm]{Problem}

\newcommand{\cF}{{\mathcal F}}
\newcommand{\cI}{{\mathcal I}}
\newcommand{\cJ}{{\mathcal J}}
\newcommand{\cS}{{\mathcal S}}
\newcommand{\cf}{{\rm cf}\/} 
\newcommand{\vare}{\varepsilon}

\begin{document}


\title[Increasing Chains]{On long increasing chains modulo flat ideals} 

\author[Saharon Shelah]{Saharon Shelah}
\email{shelah@math.huji.ac.il}

\address{The Hebrew University of Jerusalem, Einstein
Institute of Mathematics, Edmond J. Safra Campus, Givat Ram, Jerusalem
91904, Israel, and Department of Mathematics, Hill Center-Busch
Campus, Rutgers, The State University of New Jersey, 110 Frelinghuysen
Road, Piscataway, NJ 08854-8019 USA}

\thanks{I would like to thank Alice Leonhardt for the beautiful
typing. The author acknowledges support from the United States-Israel
Binational Science Foundation (Grant no. 2002323). Publication 908.}

\begin{abstract}  
  We prove that, e.g., in ${}^{(\omega_3)}(\omega_3)$ there is no sequence
  of length $\omega_4$ increasing modulo the ideal of countable sets.
\end{abstract}

\maketitle

This note is concerned with the depth of the partial order of the functions
in ${}^\kappa\gamma$ modulo the ideal of the form $\cI=[\kappa]^{<\mu}$. Let
us recall the following definitions.
\begin{dfin}
For a partial order $(P,\sqsubset)$ we define
\begin{itemize}
\item ${\rm Depth}(P,\sqsubset) = \sup\{|\cF|: {\cF}
\subseteq P$ is well--ordered by $\sqsubset\, \}$ {\em [the depth]}
\item $\cf(P,\sqsubset)=\min\{|\cF|: {\cF} \subseteq P$ is
$\sqsubset$--cofinal
which mean that for every $ p \in P $ there is
$ q \in {\mathcal F} $ such that $ p \sqsubseteq p $ 
 $\}$ {\em [the cofinality]}.

\end{itemize}
\end{dfin}
Our result (Theorem \ref{main}) states that under suitable assumptions the
depth of the partial order $({}^\kappa\gamma, <_{[\kappa]^{<\mu}})$ is at
most $|\gamma|$. In particular, letting $\mu = \aleph_1$, $\kappa = |\gamma|
= \aleph_3$ we obtain that in ${}^{(\omega_3)}(\omega_3)$ there is no
sequence of length $\omega_4$ increasing modulo the ideal of countable sets.
\medskip

Let $\kappa=\cf(\kappa)>\aleph_0$. If $\mu=\kappa$, then ${\rm
  Depth}({}^\kappa\kappa 
,<_{J^{\rm bd}_\kappa})$ can be (forced to be)
large. But for $\mu>{\rm Depth}({}^\kappa \mu, <_{J^{{\rm bd}}_\kappa})$
this implies pcf results (see \cite{Sh:410}, \cite{Sh:589}).

However, e.g., for the ideal $\cI=[\omega_3]^{\le \aleph_0}$ it is harder to
get long increasing sequence, as above for ``high $\mu$'', this leads to pcf
results 
e.g. if we assume that 
$\bar{\lambda} = \langle \lambda_i:i < \omega_3\rangle \in {}^{\omega_3}{\rm
  Reg}$, 
and in
 $(\prod\bar\lambda,<_{\cI})$ 
there is an increasing sequence moduo $ \cI $ of length say 
$>2^{\aleph_3 }+ \sup\{\lambda_i: i <\omega_3\}$ 
are much stronger than known consistency results. Even for $I =
[\omega_1]^{\le \aleph_0}$ we do not know, for $I = [\beth_\omega]^{\le
  \aleph_0}$ we know (\cite{Sh:460}), so even $[\aleph_\omega]^{\le
  \aleph_0}$ would be interesting good news. 

We hope sometime to prove, e.g.,

\begin{cject}  
For every $\mu>\theta$, in ${}^{(\theta^{+3})} \mu$ there is no increasing
sequence of length $\mu^+$ modulo $[\theta^{+3}]^{\le \theta}$.
\end{cject}

\begin{prlem}
Is it consistent
that ${}^\theta \theta$ contains $<_{\cI}$-increasing
sequence of length $\theta^+$ when $\theta =\kappa^+$ and $\cI =
[\theta]^{< \kappa}$?
\end{prlem}
\medskip

\noindent {\bf Notation:}\quad Our notation is rather standard and
compatible with that of classical textbooks (like Jech \cite{J}). 

\begin{enumerate}
\item Ordinal numbers will be denoted be the lower case initial letters of
the Greek alphabet $\alpha,\beta,\gamma,\delta,\ldots$ (with possible
subscripts). Cardinal numbers will be called $\kappa,\lambda,\mu,\theta$.  
\item For a set $X$ and a cardinal $\theta$, $[X]^{\theta}$ (or
$[X]^{<\theta}$, respectively) stands for the family of subsets of $X$
of size $\theta$ ($<\theta$, respectively). 
\end{enumerate}

\begin{thrm}
\label{main}
Assume $\mu^+< \kappa \le \theta$ and $\cJ=[\kappa]^{<\mu}$ and 
$\cf([\theta]^\mu,\subseteq) \le \theta$. Let $\gamma<\theta^+$. 
Then ${\rm Depth}({}^{\kappa}\gamma,<_{\cJ})\leq \theta$, i.e.,
there is no $<_{\cJ}$--increasing sequence $\langle f_\alpha:\alpha <
\theta^+\rangle$ of functions from ${}^\kappa\gamma$ modulo $\cJ$.
\end{thrm}

\begin{proof} 
Assume towards contradiction that there is a $<_\cJ$--increasing sequence
$\langle f_ \zeta : \zeta <\theta^+\rangle \subseteq {}^\kappa\gamma$.

Let $\cS\subseteq [\gamma]^\mu$ be cofinal of cardinality $\le \theta$
(exists as $|\gamma|\leq\theta$ and $\cf([\theta]^\mu,\subseteq) \le
\theta$).  For every $s \in {\cS}$ and $\beta<\kappa$, $\zeta<\theta^+$
we let  
\begin{itemize}
\item $I(\beta) := [\beta,\beta + \mu)$,
\item $f^s_\zeta \in {}^\kappa(\gamma+1)$ be defined by
$f^s_\zeta(i) = \min(s \cup \{\gamma\} \setminus f_\zeta(i))$, 
\item $f^{s,\beta}_\zeta \in {}^{I(\beta)}(\gamma+1)$ be defined as
$f^s_\zeta \restriction I(\beta)$.
\end{itemize}
Now, for each $s\in \cS$ we have
\begin{enumerate}
\item[$(*)_1$]  
\begin{enumerate}
\item[(a)]  for every $\zeta<\theta^+$, $f^{s,\beta}_\zeta:{I(\beta)}
  \longrightarrow s\cup\{\gamma\}$, 
\item[(b)]  if $\zeta<\xi<\theta^+$, then $f^{s,\beta}_\zeta \le
  f^{s,\beta}_\xi \mod [I(\beta)]^{<\mu}$. 
\end{enumerate}
\end{enumerate}
For $s \in\cS$ we define
\begin{enumerate}
\item[$(*)_2$] $B_s=\{\beta < \kappa:(\forall\zeta<\theta^+)(\exists \xi>
  \zeta) \neg(f^{s,\beta}_\zeta = f^{s,\beta}_\xi \mod [I(\beta)]^{<\mu})\,\}$.
\end{enumerate}
Plainly, we may choose a sequence $\langle C^s_\beta:\beta<\kappa,\ s\in\cS
\rangle$ such that
\begin{enumerate}
\item[$(*)_3$]  
\begin{enumerate}
\item[(a)] $C^s_\beta$ is a club of $\theta^+$,
\item[(b)] if $\beta\in B_s$ and $\xi,\zeta\in C^s_\beta$ are such that
  $\zeta < \xi$, then $\neg(f^{s,\beta}_\zeta=f^{s,\beta}_\xi \mod
  [I(\beta)]^{<\mu})$, 
\item[(c)] if $\beta\in\kappa\setminus B_s$, then $f^{s,\beta}_\zeta =
  f^{s,\beta}_\xi \mod   [I(\beta)]^{<\mu}$ whenever
  $\min(C^s_\beta)\le\zeta\le\xi <\theta^+$.
\end{enumerate}
\end{enumerate}
Then, as $|\cS|\le \theta$ and $\kappa\le\theta$, we have
\begin{enumerate}
\item[$(*)_4$]  the set $C:= \bigcap\{C^s_\beta:s \in\cS$ and
  $\beta<\kappa\}$ is a club of $\theta^+$ .  
\end{enumerate}
Choose a sequence $\langle\alpha_\vare:\vare<\mu^+\rangle \subseteq C$
increasing with $\varepsilon$. Then, for all $\varepsilon<\zeta<\mu^+$,  
\begin{enumerate}
\item[$(*)_5$] $u_{\vare,\zeta}:=\{i<\kappa: f_{\alpha_\vare}(i)\ge
  f_{\alpha_\zeta} (i)\}\in\cJ$.
\end{enumerate}
We have assumed that $\mu^+<\kappa$, so we can find $\delta<\kappa$ such
that  
\begin{enumerate}
\item[$(*)_6$]
\begin{enumerate}
\item[(a)] $I(\delta)= [\delta,\delta+\mu)$ is disjoint from
  $\bigcup\{u_{\vare,\zeta}:\vare<\zeta<\mu^+\}$, and hence
\item[(b)] the sequence $\langle f_{\alpha_\vare}(i):\vare<\mu^+\rangle$ is 
increasing for each $i \in I(\delta)$.
\end{enumerate}
\end{enumerate}
As $|I(\delta)|=\mu$ and $\cS\subseteq [\gamma]^{\le\mu}$ is cofinal (for
the partial order $\subseteq$), we can find $s\in\cS$ such that
\begin{enumerate}
\item[$(*)_7$] $\{f_{\alpha_0}(i),f_{\alpha_1}(i):i\in I(\delta)\}\subseteq
  s$. 
\end{enumerate}
It follows from $(*)_6 + (*)_7$ that for every $i \in I(\delta)$
\begin{enumerate}
\item[$(*)_8$]  $f^s_{\alpha_0}(i)=f_{\alpha_0}(i)<f_{\alpha_1}(i) =
  f^s_{\alpha_1}(i)$. 
\end{enumerate}
As $\alpha_0<\alpha_1$ are from $C$ and $I(\delta)\notin \cJ$, recalling
$(*)_2 + (*)_3 + (*)_4$, clearly
\begin{enumerate}
\item[$(*)_9$]  $\delta\in B_s$.
\end{enumerate}
Therefore, as $\alpha_\vare\in C\subseteq C^s_\delta$ for $\vare<\mu^+$ and
$\alpha_\vare$ is increasing with $\vare$, we have
\begin{enumerate}
\item[$(*)_{10}$]  for every $\vare<\mu^+$ there is $i_\vare\in I(\delta)$
  such that 
\begin{enumerate}
\item[$(\alpha)$]  $f^s_{\alpha_\vare}(i_\vare)<f^s_{\alpha_{\vare+1}}
  (i_\vare)$,
\end{enumerate}
and hence there is $j_\vare\in s$ such that
\begin{enumerate}
\item[$(\beta)$] $f^s_{\alpha_\vare}(i_\vare) \le j_\vare<
  f^s_{\alpha_{\vare +1}}(i_\vare)$ 
\end{enumerate}
and therefore 
\begin{enumerate}
\item[$(\gamma)$] $f_{\alpha_\vare}(i_\vare) \le j_\vare < f_{\alpha_{\vare
      +1}}(i_\vare)$. 
\end{enumerate}
\end{enumerate}
But $|I(\delta)| + |s| = \mu < \mu^+$, so for some pair $(j_*,i_*) \in s
\times I(\delta)$ we may choose $\vare_1< \vare_2<\mu^+$ such that  
\begin{enumerate}
\item[$(*)_{11}$]  $j_{\vare_1} =j_{\vare_2} = j_*$ and $i_{\vare_1}
  =i_{\vare_2} = i_*$. 
\end{enumerate}
But the sequence $\langle f_{\alpha_\vare}(i_*):\vare<\theta^+\rangle$ is
increasing by $(*)_6(b)$ (see the choice of $\delta$), so 
\[f_{\alpha_{\vare_1}}(i_*)<f_{\alpha_{\vare_1+1}}(i_*) \le
f_{\alpha_{\vare_2}}(i_*)< f_{\alpha_{\vare_2+1}}(i_*).\]
It follows from $(*)_{10}(\gamma)+(*)_{11}$ that the ordinal $j_*$ belongs
to $[f_{\alpha_{\vare_1}}(i_*),f_{\alpha_{\vare_1+1}}(i_*))$ and to
$[f_{\alpha_{\vare_2}}(i_*),f_{\alpha_{\vare_2+1}}(i_*))$, which are
disjoint intervals, a contradiction.  
\end{proof}

Similarly,

\begin{thrm}  
Assume that
\begin{enumerate}
\item[(a)] $\cJ$ is an ideal on $\kappa$,
\item[(b)] $I_\beta \in [\kappa]^\mu$, $I_\beta \notin\cJ$ for
  $\beta<\kappa$, 
\item[(c)] $\theta = |\gamma|+\kappa$ and
  $\cf([\theta]^\mu,\subseteq)<\lambda$,
\item[(d)] if $u_\vare\in\cJ$ for $\vare<\mu^+$, then for some $\beta <
\kappa$ the set $I_\beta$ is disjoint from $\bigcup_{\vare<\mu^+}
u_\varepsilon$. 
\end{enumerate}
Then there is no $<_{\cJ}$--increasing sequence of functions from $\kappa$
to $\gamma$ of length $\lambda$. 
\end{thrm}

\begin{proof}  
Without loss of generality $\lambda$ is the successor of
$\cf([\theta]^\mu,\subseteq)$ hence is regular.  The proof is similar to the 
proof of Theorem \ref{main}. 
\end{proof}

\end{document}